# AN INTRODUCTION TO THE NEUTROSOPHIC PROBABILITY APPLIED IN QUANTUM PHYSICS

Florentin Smarandache, Ph. D.
University of New Mexico
Gallup, NM 87301, USA

**Abstract.**
In this paper one generalizes the *classical probability* and *imprecise probability* to the notion of "neutrosophic probability" in order to be able to model Heisenberg's Uncertainty Principle of a particle's behavior, Schrödinger's Cat Theory, and the state of bosons which do not obey Pauli's Exclusion Principle (in quantum physics). Neutrosophic probability is close related to neutrosophic logic and neutrosophic set, and etymologically derived from "neutrosophy" [58, 59].

## 1. Introduction.

One consequence of the Heisenberg's Uncertainty Principle says that it is impossible to fully predict the behavior of a particle, also the causality principle cannot apply at the atomic level.

For example the Schrödinger's Cat Theory says that the quantum state of a photon can basically be in more than one place in the same time which, translated to the neutrosophic set, means that an element (quantum state) belongs and does not belong to a set (a place) in the same time; or an element (quantum state) belongs to two different sets (two different places) in the same time. It is a question of "alternative worlds" theory very well represented by the neutrosophic set theory.

In Schrödinger's Equation on the behavior of electromagnetic waves and "matter waves" in quantum theory, the wave function < which describes the superposition of possible states may be simulated by a neutrosophic function, i.e. a function whose values are not unique for each argument from the domain of definition (the vertical line test fails, intersecting the graph in more points).

How to describe a particle $\zeta$ in the infinite micro-universe that belongs to two distinct places $P_1$ and $P_2$ in the same time? $\zeta \in P_1$ and $\zeta \notin P_1$ as a true contradiction, or $\zeta \in P_1$ and $\zeta \in \neg P_1$.

Or, how to describe two distinct bosons $b_1$ and $b_2$, which do not obey Pauli's Exclusion Principle, i.e. they belong to the same quantum or energy state in the same time?

Or, how to calculate the truth-value of Zen (in Japanese) / Chan (in Chinese) doctrine philosophical proposition: the present is eternal and comprises in itself the past and the future?

In Eastern Philosophy the contradictory utterances form the core of the Taoism and Zen/Chan (which emerged from Buddhism and Taoism) doctrines.
How to judge the truth-value of a metaphor, or of an ambiguous statement, or of a social phenomenon which is positive from a standpoint and negative from another standpoint?

We better describe them, using the attribute "neutrosophic" than "fuzzy" or any other, a quantum particle that neither exists nor non-exists.

## 2. Non-Standard Real Numbers and Non-Standard Real Sets.

Let T, I, F be standard or non-standard real subsets $\subseteq\ ]^-0, 1^+[$,
with   sup T = t_sup,  inf T = t_inf,
       sup I = i_sup,  inf I = i_inf,
       sup F = f_sup,  inf F = f_inf,
and n_sup = t_sup + i_sup + f_sup,
    n_inf = t_inf + i_inf + f_inf.
Obviously: t_sup, i_sup, f_sup $\leq 1^+$, and t_inf, i_inf, f_inf $\geq\ ^-0$,
whereas n_sup $\leq 3^+$ and n_inf $\geq\ ^-0$.

The subsets T, I, F are not necessarily intervals, but may be any real subsets: discrete or continuous; single-element, finite, or (either countably or uncountably) infinite; union or intersection of various subsets; etc.
They may also overlap. These real subsets could represent the relative errors in determining t, i, f (in the case when the subsets T, I, F are reduced to points).

This representation is closer to the human mind reasoning. It characterizes/catches the *imprecision* of knowledge or linguistic inexactitude received by various observers (that's why T, I, F are subsets - not necessarily single-elements), *uncertainty* due to incomplete knowledge or acquisition errors or stochasticity (that's why the subset I exists), and *vagueness* due to lack of clear contours or limits (that's why T, I, F are subsets and I exists; in particular for the appurtenance to the neutrosophic sets).
One has to specify the superior (x_sup) and inferior (x_inf) limits of the subsets because in many problems arises the necessity to compute them.

The real number x is said to be infinitesimal if and only if for all positive integers n one has $|x| < 1/n$. Let $\varepsilon > 0$ be a such infinitesimal number. The *hyper-real number set* is an extension of the real number set, which includes classes of infinite numbers and classes of infinitesimal numbers. Let's consider the non-standard finite numbers $1^+ = 1+\varepsilon$, where "1" is its standard part and "$\varepsilon$" its non-standard part, and $^-0 = 0-\varepsilon$, where "0" is its standard part and "$\varepsilon$" its non-standard part.
Then, we call $]^-0, 1^+[$ a non-standard unit interval. Obviously, 0 and 1, and analogously non-standard numbers infinitely small but less than 0 or infinitely small but greater than 1, belong to the non-standard unit interval. Actually, by "$^-a$" one signifies a monad, i.e. a set of hyper-real numbers in non-standard analysis:
   . $(^-a) = \{a-x: x \in \mathbb{R}^*,\ x\ \text{is infinitesimal}\}$,
and similarly "$b^+$" is a monad:
   . $(b^+) = \{b+x: x \in \mathbb{R}^*,\ x\ \text{is infinitesimal}\}$.

Generally, the left and right borders of a non-standard interval $](^-a, b^+[$ are vague, imprecise, themselves being non-standard (sub)sets $\mu(^-a)$ and $\mu(b^+)$ as defined above.

Combining the two before mentioned definitions one gets, what we would call, a binad of "$^-c^+$":

$\mu(^-c^+) = \{c-x: x \in \Re^*, x \text{ is infinitesimal}\} \cup \{c+x: x \in \Re^*, x \text{ is infinitesimal}\}$, which is a collection of open punctured neighborhoods (balls) of c.

Of course, $^-a < a$ and $b^+ > b$. No order between $^-c^+$ and c.

Addition of non-standard finite numbers with themselves or with real numbers:

$^-a + b = ^-(a + b)$
$a + b^+ = (a + b)^+$
$^-a + b^+ = ^-(a + b)^+$
$^-a + ^-b = ^-(a + b)$ (the left monads absorb themselves)
$a^+ + b^+ = (a + b)^+$ (analogously, the right monads absorb themselves)

Similarly for subtraction, multiplication, division, roots, and powers of non-standard finite numbers with themselves or with real numbers.

By extension let inf $](^-a, b^+[ = ^-a$ and sup $](^-a, b^+[ = b^+$.

### 3. A Logical Connection.

Łukasiewicz, together with Kotarbiński and Leśniewski from the Warsaw Polish Logic group (1919-1939), questioned the status of truth: eternal, sempiternal (everlasting, perpetual), or both?

Let's borrow from the modal logic the notion of "world", which is a semantic device of what the world might have been like. Then, one says that the neutrosophic truth-value of a statement A, $NL_t(A) = 1^+$ if A is 'true in all possible worlds' (syntagme first used by Leibniz) and all conjunctures, that one may call "absolute truth" (in the modal logic it was named *necessary truth*, Dinulescu-Câmpina [9] names it 'intangible absolute truth' ), whereas $NL_t(A) = 1$ if A is true in at least one world at some conjuncture, we call this "relative truth" because it is related to a 'specific' world and a specific conjuncture (in the modal logic it was named *possible truth*). Because each 'world' is dynamic, depending on an ensemble of parameters, we introduce the sub-category 'conjuncture' within it to reflect a particular state of the world.

How can we differentiate <the truth behind the truth>? What about the <metaphoric truth>, which frequently occurs in the humanistic field? Let's take the proposition "99% of the politicians are crooked" (Sonnabend [60], Problem 29, p. 25). "No," somebody furiously comments, "100% of the politicians are crooked, *even more*!" How do we interpret this "even more" (than 100%), i. e. more than the truth?

One attempts to formalize. For $n \geq 1$ one defines the "n-level relative truth" of the statement A if the statement is true in at least n distinct worlds, and similarly "countably-" or "uncountably-level relative truth" as gradual degrees between "first-level relative truth" (1) and "absolute truth" ($1^+$) in the monad $\mu(1^+)$. Analogue definitions one gets by substituting "truth" with "falsehood" or "indeterminacy" in the above.

In *largo sensu* the notion "world" depends on parameters, such as: space, time, continuity, movement, modality, (meta)language levels, interpretation, abstraction, (higher-order) quantification, predication, complement constructions, subjectivity, context, circumstances, etc. Pierre d'Ailly upholds that the truth-value of a proposition depends on the sense, on the metaphysical level, on the language and meta-language; the auto-reflexive propositions (with reflection on themselves) depend on the mode of representation (objective/subjective, formal/informal, real/mental).

In a formal way, let's consider the world W as being generated by the formal system FS. One says that statement A belongs to the world W if A is a well-formed formula (*wff*) in W, i.e. a string of symbols from the alphabet of W that conforms to the grammar of the formal language endowing W. The grammar is conceived as a set of functions (formation rules) whose inputs are symbols strings and outputs "yes" or "no". A formal system comprises a formal language (alphabet and grammar) and a deductive apparatus (axioms and/or rules of inference). In a formal system the rules of inference are syntactically and typographically formal in nature, without reference to the meaning of the strings they manipulate.

Similarly for the neutrosophic falsehood-value, $NL_f(A) = 1^+$ if the statement A is false in all possible worlds, we call it "absolute falsehood", whereas $NL_f(A) = 1$ if the statement A is false in at least one world, we call it "relative falsehood". Also, the neutrosophic indeterminacy-value $NL_i(A) = 1^+$ if the statement A is indeterminate in all possible worlds, we call it "absolute indeterminacy", whereas $NL_i(A) = 1$ if the statement A is indeterminate in at least one world, we call it "relative indeterminacy".

On the other hand, $NL_t(A) = {}^-0$ if A is false in all possible world, whereas $NL_t(A) = 0$ if A is false in at least one world; $NL_f(A) = {}^-0$ if A is true in all possible world, whereas $NL_f(A) = 0$ if A is true in at least one world; and $NL_i(A) = {}^-0$ if A is indeterminate in no possible world, whereas $NL_i(A) = 0$ if A is not indeterminate in at least one world.

The $^-0$ and $1^+$ monads leave room for degrees of super-truth (truth whose values are greater than 1), super-falsehood, and super-indeterminacy.

Here there are some corner cases:

There are tautologies, some of the form "B is B", for which $NL(B) = (1^+, {}^-0, {}^-0)$, and contradictions, some of the form "C is not C", for which $NL(B) = ({}^-0, {}^-0, 1^+)$.

While for a paradox, P, $NL(P) = (1,1,1)$. Let's take the Epimenides Paradox, also called the Liar Paradox, "This very statement is true". If it is true then it is false, and if it is false then it is true. But the previous reasoning, due to the contradictory results, indicates a high indeterminacy too. The paradox is the only proposition true and false in the same time in the same world, and indeterminate as well!

Let's take the Grelling's Paradox, also called the heterological paradox [Suber, 1999], "If an adjective truly describes itself, call it 'autological', otherwise call it 'heterological'. Is 'heterological' heterological?" Similarly, if it is, then it is not; and if it is not, then it is.

For a not well-formed formula, nwff, i.e. a string of symbols which do not conform to the syntax of the given logic, $NL(nwff) = n/a$ (undefined). A proposition which may not be considered a proposition was called by the logician Paulus Venetus *flatus voci*. $NL(flatus\ voci) = n/a$.

## 4. Operations with Standard and Non-Standard Real Subsets.

Let $S_1$ and $S_2$ be two (unidimensional) standard or non-standard real subsets, then one defines:

Addition of sets:
$S_1 \oplus S_2 = \{x \mid x=s_1+s_2,\ \text{where}\ s_1 \in S_1\ \text{and}\ s_2 \in S_2\}$,
with $\inf S_1 \oplus S_2 = \inf S_1 + \inf S_2$, $\sup S_1 \oplus S_2 = \sup S_1 + \sup S_2$;
and, as some particular cases, we have
$\{a\} \oplus S_2 = \{x \mid x=a+s_2,\ \text{where}\ s_2 \in S_2\}$

with inf $\{a\}\oplus S_2 = a + \inf S_2$, sup $\{a\}\oplus S_2 = a + \sup S_2$;
also $\{1\}\oplus S_2 = \{x \mid x=1+s_2, \text{ where } s_2\in S_2\}$
with inf $\{1\}\oplus S_2 = 1 + \inf S_2$, sup $\{1\}\oplus S_2 = 1 + \sup S_2$.

Subtraction of sets:
$S_1\ominus S_2 = \{x \mid x=s_1-s_2, \text{ where } s_1\in S_1 \text{ and } s_2\in S_2\}$.
For real positive subsets (most of the cases will fall in this range) one gets
inf $S_1\ominus S_2 = \inf S_1 - \sup S_2$, sup $S_1\ominus S_2 = \sup S_1 - \inf S_2$;
and, as some particular cases, we have
$\{a\}\ominus S_2 = \{x \mid x=a-s_2, \text{ where } s_2\in S_2\}$,
with inf $\{a\}\ominus S_2 = a - \sup S_2$, sup $\{a\}\ominus S_2 = a - \inf S_2$;
also $\{1\}\ominus S_2 = \{x \mid x=1-s_2, \text{ where } s_2\in S_2\}$,
with inf $\{1\}\ominus S_2 = 1 - \sup S_2$, sup $\{1\}\ominus S_2 = 100 - \inf S_2$.

Multiplication of sets:
$S_1\odot S_2 = \{x \mid x=s_1\cdot s_2, \text{ where } s_1\in S_1 \text{ and } s_2\in S_2\}$.
For real positive subsets (most of the cases will fall in this range) one gets
inf $S_1\odot S_2 = \inf S_1 \cdot \inf S_2$, sup $S_1\odot S_2 = \sup S_1 \cdot \sup S_2$;
and, as some particular cases, we have
$\{a\}\odot S_2 = \{x \mid x=a\cdot s_2, \text{ where } s_2\in S_2\}$,
with inf $\{a\}\odot S_2 = a * \inf S_2$, sup $\{a\}\odot S_2 = a \cdot \sup S_2$;
also $\{1\}\odot S_2 = \{x \mid x=1\cdot s_2, \text{ where } s_2\in S_2\}$,
with inf $\{1\}\odot S_2 = 1 \cdot \inf S_2$, sup $\{1\}\odot S_2 = 1 \cdot \sup S_2$.

Division of a set by a number:
Let $k \in \mathbb{R}^*$, then $S_1\oslash k = \{x \mid x=s_1/k, \text{ where } s_1\in S_1\}$,

Let $(T_1, I_1, F_1)$ and $(T_2, I_2, F_2)$ be standard or non-standard triplets of real subsets which $\in P(\mathbb{K}\ ^-0, 1^+\mathbb{M})^3$, where $P(\mathbb{K}\ ^-0, 1^+\mathbb{M})$ is the set of all subsets of non-standard unit interval $\mathbb{K}\ ^-0, 1^+\mathbb{M}$, then we define:
$(T_1, I_1, F_1) \boxplus (T_2, I_2, F_2) = (T_1\oplus T_2, I_1\oplus I_2, F_1\oplus F_2)$,
$(T_1, I_1, F_1) \boxminus (T_2, I_2, F_2) = (T_1\ominus T_2, I_1\ominus I_2, F_1\ominus F_2)$,
$(T_1, I_1, F_1) \boxdot (T_2, I_2, F_2) = (T_1\odot T_2, I_1\odot I_2, F_1\odot F_2)$.

## 5. Neutrosophic Probability:
Is a generalization of the classical probability in which the chance that an event A occurs is t% true - where t varies in the subset T, i% indeterminate - where i varies in the subset I, and f% false - where f varies in the subset F.
One notes $NP(A) = (T, I, F)$.
It is also a generalization of the imprecise probability, which is an interval-valued distribution function.

## 6. Neutrosophic Statistics:
Is the analysis of the events described by the neutrosophic probability.
This is also a generalization of the classical statistics and imprecise statistics.

## 7. Neutrosophic Probability Space.

The universal set, endowed with a neutrosophic probability defined for each of its subset, forms a neutrosophic probability space.

Let A and B be two neutrosophic events, and $NP(A) = (T_1, I_1, F_1)$, $NP(B) = (T_2, I_2, F_2)$ their neutrosophic probabilities. Then we define:

$NP(A \cap B) = NP(A) \boxdot NP(B)$.
$NP(\neg A) = \{1\} \boxminus NP(A)$.
$NP(A \cup B) = NP(A) \boxplus NP(B) \boxminus NP(A) \boxdot NP(B)$.

1. $NP(\text{impossible event}) = (T_{imp}, I_{imp}, F_{imp})$,
where sup $T_{imp} \leq 0$, inf $F_{imp} \geq 1$; no restriction on $I_{imp}$.
   $NP(\text{sure event}) = (T_{sur}, I_{sur}, F_{sur})$,
where inf $T_{sur} \geq 1$, sup $F_{sur} \leq 0$; no restriction on $I_{sur}$.
   $NP(\text{totally indeterminate event}) = (T_{ind}, I_{ind}, F_{ind})$;
where inf $I_{ind} \geq 1$; no restrictions on $T_{ind}$ or $F_{ind}$.
2. $NP(A) \in \{(T, I, F)$, where T, I, F are real subsets which may overlap$\}$.
3. $NP(A \cup B) = NP(A) \boxplus NP(B) \boxminus NP(A \cap B)$.
4. $NP(A) = \{1\} \boxminus NP(\neg A)$.

## 8. Applications:

#1. From a pool of refugees, waiting in a political refugee camp in Turkey to get the American visa, a% have the chance to be accepted - where a varies in the set A, r% to be rejected - where r varies in the set R, and p% to be in pending (not yet decided) - where p varies in P.
Say, for example, that the chance of someone Popescu in the pool to emigrate to USA is (between) 40-60% (considering different criteria of emigration one gets different percentages, we have to take care of all of them), the chance of being rejected is 20-25% or 30-35%, and the chance of being in pending is 10% or 20% or 30%. Then the neutrosophic probability that Popescu emigrates to the Unites States is
   NP(Popescu) = ( (40-60), (20-25)U(30-35), {10,20,30} ), closer to the life's thinking.
This is a better approach than the classical probability, where $40 \leq P(\text{Popescu}) \leq 60$, because from the pending chance - which will be converted to acceptance or rejection - Popescu might get extra percentage in his will to emigration,
and also the superior limit of the subsets sum
   60+35+30 > 100
and in other cases one may have the inferior sum < 0,
while in the classical fuzzy set theory the superior sum should be 100 and the inferior sum $\geq 0$.
In a similar way, we could say about the element Popescu that
Popescu( (40-60), (20-25)U(30-35), {10,20,30} ) belongs to the set of accepted refugees.
#2. The probability that candidate C will win an election is say 25-30% true (percent of people voting for him), 35% false (percent of people voting against him), and 40% or 41% indeterminate (percent of people not coming to the ballot box, or giving a blank vote - not selecting anyone, or giving a negative vote - cutting all candidates on the list).
Dialectic and dualism don't work in this case anymore.

#3. Another example, the probability that tomorrow it will rain is say 50-54% true according to meteorologists who have investigated the past years' weather, 30 or 34-35% false according to today's very sunny and droughty summer, and 10 or 20% undecided (indeterminate).

#4. The probability that Yankees will win tomorrow versus Cowboys is 60% true (according to their confrontation's history giving Yankees' satisfaction), 30-32% false (supposing Cowboys are actually up to the mark, while Yankees are declining), and 10 or 11 or 12% indeterminate (left to the hazard: sickness of players, referee's mistakes, atmospheric conditions during the game).  These parameters act on players' psychology.

## 9. Remarks:

Neutrosophic probability is useful to those events which involve some degree of indeterminacy (unknown) and more criteria of evaluation - as quantum physics.  This kind of probability is necessary because it provides a better representation than classical probability to uncertain events.

## 10. Generalizations of Other Probabilities.

In the case when the truth- and falsity-components are complementary, i.e. no indeterminacy and their sum is 1, one falls to the classical probability.  As, for example, tossing dice or coins, or drawing cards from a well-shuffled deck, or drawing balls from an urn.

An interesting particular case is for n=1, with $0 \leq t,i,f \leq 1$, which is closer to the classical probability.

For n=1 and i=0, with $0 \leq t,f \leq 1$, one obtains the classical probability.

From the intuitionistic logic, paraconsistent logic, dialetheism, faillibilism, paradoxism, pseudoparadoxism, and tautologism we transfer the  "adjectives" to probabilities, i.e. we define the intuitionistic probability (when the probability space is incomplete), paraconsistent probability, faillibilist probability, dialetheist probability, paradoxist probability, pseudoparadoxist probability, and tautologic probability respectively.

Hence, the neutrosophic probability generalizes:
- the intuitionistic probability, which supports incomplete (not completely known/determined) probability spaces (for 0<n<1 and i=0, $0 \leq t,f \leq 1$) or incomplete events whose probability we need to calculate;
- the classical probability (for n=1 and i=0, and $0 \leq t,f \leq 1$);
- the paraconsistent probability (for n>1 and i=0, with both t,f<1);
- the dialetheist probability, which says that intersection of some disjoint probability spaces is not empty (for t=f=1 and i=0; some paradoxist probabilities can be denoted this way);
- the faillibilist probability (for i>0);
- the pseudoparadoxism (for n_sup>1 or n_inf<0);
- the tautologism (for t_sup>1).

Compared with all other types of classical probabilities, the neutrosophic probability introduces a percentage of "indeterminacy" - due to unexpected parameters hidden in some probability spaces, and let each component t, i, f be even boiling *over* 1 to $1^+$ (overflooded) or freezing *under* 0 (underdried) to $^-0$.

For example: an element in some tautological probability space may have t>1, called "overprobable" (i.e. t = $1^+$).  Similarly, an element in some paradoxist probability space may be "overindeterminate" (for i>1), or "overunprobable" (for f>1, in some unconditionally false appurtenances); or "underprobable" (for t<0, i.e. t = ⁻0, in some unconditionally false appurtenances), "underindeterminate" (for i<0, in some unconditionally true or false appurtenances), "underunprobable" (for f<0, in some unconditionally true appurtenances).

This is because we should make a distinction between unconditionally true (t>1, and f<0 or i<0) and conditionally true appurtenances (t≤1, and f≤1 or i≤1).

## 11. Other Examples.

Let's consider a neutrosophic set a collection of possible locations (positions) of particle x. And let A and B be two neutrosophic sets.

One can say, by language abuse, that any particle x neutrosophically belongs to any set, due to the percentages of truth/indeterminacy/falsity involved, which varies between $^-0$ and $1^+$. For example: x(0.5, 0.2, 0.3) belongs to A (which means, with a probability of 50% particle x is in a position of A, with a probability of 30% x is not in A, and the rest is undecidable); or y(0, 0, 1) belongs to A (which normally means y is not for sure in A); or z(0, 1, 0) belongs to A (which means one does know absolutely nothing about z's affiliation with A).

More general, x( (0.2-0.3), (0.40-0.45)∪[0.50-0.51], {0.2, 0.24, 0.28} ) belongs to the set A, which means:

- with a probability in between 20-30% particle x is in a position of A (one cannot find an exact approximate because of various sources used);
- with a probability of 20% or 24% or 28% x is not in A;
- the indeterminacy related to the appurtenance of x to A is in between 40-45% or between 50-51% (limits included).

The subsets representing the appurtenance, indeterminacy, and falsity may overlap, and n_sup = 30%+51%+28% > 100% in this case.